\documentclass[english,twoside]{article}

\usepackage[german,english]{babel}
\usepackage{amsmath}
\usepackage{amssymb,latexsym}
\usepackage{longtable}
\usepackage{booktabs}
\usepackage{tabularx}
\usepackage{multicol}
\usepackage{url}

\newtheorem{deff}{Definition}

\newtheorem{thm}[deff]{Theorem}
\newtheorem{cor}[deff]{Corollary}

\newtheorem{probl}[deff]{Problem}

\makeatletter

\title{Combinatorial $3$-Manifolds with $10$ Vertices}

\author{\Large Frank H.~Lutz}

\date{}

\begin{document}

\selectlanguage{english}

\maketitle

\bigskip

\begin{abstract}
\mbox{}
We give a complete enumeration of all combinatorial $3$-manifolds with $10$~vertices:
There are precisely $247882$ triangulated $3$-spheres with $10$ vertices
as well as $518$ vertex-minimal triangulations of the sphere product $S^2\!\times\!S^1$
and $615$ triangulations of the twisted sphere product $S^2\hbox{$\times\hspace{-1.62ex}\_\hspace{-.4ex}\_\hspace{.7ex}$}S^1$.

All the $3$-spheres with up to $10$ vertices are shellable, 
but there are $29$ vertex-minimal non-shellable $3$-balls with $9$ vertices.
\mbox{}
\end{abstract}

\section{Introduction}

Let $M$ be a triangulated $3$-manifold with $n$ vertices and face vector $f=(n,f_1,f_2,f_3)$. 
By Euler's equation, $n - f_1 + f_2 - f_3 = 0$, and by double counting
the edges of the ridge-facet incidence graph, $2 f_2 = 4 f_3$,
it follows that
\begin{equation}
f=(n,f_1,2f_1-2n,f_1-n).
\end{equation}

A complete characterization of the $f$-vectors of
the $3$-sphere $S^3$, the sphere product $S^2\!\times\!S^1$,
the twisted sphere product (or $3$-dimensional Klein bottle)
$S^2\hbox{$\times\hspace{-1.62ex}\_\hspace{-.4ex}\_\hspace{.7ex}$}S^1$,
and of the real projective $3$-space ${\mathbb R}{\bf P}^{\,3}$
was given by Walkup.
\begin{thm} {\rm (Walkup \cite{Walkup1970})}\label{thm:Walkup_gamma}
For every $3$-manifold $M$ there is an integer $\gamma (M)$
such that
\begin{equation}
f_1\geq 4n+\gamma (M)
\end{equation}
for every triangulation of $M$ with $n$ vertices and $f_1$ edges. 
Moreover, there is an integer $\gamma^*(M)\geq \gamma (M)$
such that for every pair $(n,f_1)$ with $n\geq 5$ and
\begin{equation}
\binom{n}{2}\geq f_1\geq 4n+\gamma^*(M)
\end{equation}
there is a triangulation of $M$ with $n$ vertices and $f_1$ edges.
In particular,
\begin{itemize}
\item[{\rm (a)}] $\gamma^*=\gamma =-10$\, for\, $S^3$,
\item[{\rm (b)}] $\gamma^*=\gamma =0$\, for\, $S^2\hbox{$\times\hspace{-1.62ex}\_\hspace{-.4ex}\_\hspace{.7ex}$}S^1$,
\item[{\rm (c)}] $\gamma^*=1$\, and\, $\gamma =0$\, for\,
  $S^2\!\times\!S^1$, where, with the exception $(9,36)$, all pairs $(n,f_1)$
  with $n\geq 5$ and $\binom{n}{2}\geq f_1\geq 4n$ occur,
\item[{\rm (d)}] $\gamma^*=\gamma =7$\, for\, ${\mathbb R}{\bf P}^{\,3}$, and
\item[{\rm (e)}] $\gamma^*(M)\geq\gamma (M)\geq 8$\, for all
  other $3$-manifolds $M$.
\end{itemize}
\end{thm}
By Walkup's theorem, vertex-minimal triangulations
of $S^2\hbox{$\times\hspace{-1.62ex}\_\hspace{-.4ex}\_\hspace{.7ex}$}S^1$,
$S^2\!\times\!S^1$, and ${\mathbb R}{\bf P}^{\,3}$ have $9$, $10$, and
$11$ vertices, respectively. The $3$-sphere can be triangulated
vertex-minimally as the boundary of the $4$-simplex with $5$ vertices.
But otherwise, rather little is known on vertex-minimal
triangulations of $3$- and higher-dimensional manifolds.
See \cite{Lutz2003bpre}, \cite{Lutz2005bpre}, and \cite{Swartz2005pre}
for a discussion, further references, and for various examples of small triangulations
of $3$-mani\-folds.

The exact numbers of different combinatorial types of triangulations
of $S^3$, $S^2\hbox{$\times\hspace{-1.62ex}\_\hspace{-.4ex}\_\hspace{.7ex}$}S^1$,
and $S^2\!\times\!S^1$ with up to $9$ vertices and of neighborly
triangulations (i.e., triangulations with complete $1$-skeleton)
with $10$ vertices were obtained by
\begin{tabbing}
      oooooooooooooooopppppppppppp\hspace{1mm} \=  \kill
      Gr\"unbaum and Sreedharan \cite{GruenbaumSreedharan1967} \> (simplicial $4$-polytopes with $8$ vertices),\\
      Barnette \cite{Barnette1973c}                            \> (combinatorial $3$-spheres with $8$ vertices),\\
      Altshuler \cite{Altshuler1974}                           \> (combinatorial $3$-manifolds with up to $8$ vertices),\\
      Altshuler and Steinberg \cite{AltshulerSteinberg1973}    \> (neighborly $4$-polytopes with $9$ vertices),\\
      Altshuler and Steinberg \cite{AltshulerSteinberg1974}    \> (neighborly $3$-manifolds with $9$ vertices),\\
      Altshuler and Steinberg \cite{AltshulerSteinberg1976}    \> (combinatorial $3$-manifolds with $9$ vertices),\\
      Altshuler \cite{Altshuler1977}                           \> (neighborly $3$-manifolds with $10$ vertices).
\end{tabbing}

\begin{table}
\small\centering
\defaultaddspace=0.15em
\caption{Combinatorial $3$-manifolds with up to $10$ vertices.}\label{tbl:ten3d_leq10}
\begin{tabular*}{10cm}{@{\extracolsep{\fill}}c@{\hspace{5mm}}r@{\hspace{5mm}}r@{\hspace{5mm}}r@{\hspace{5mm}}r@{}}
 \addlinespace
 \addlinespace
 \addlinespace
 \addlinespace
 \addlinespace
 \addlinespace
\toprule
 \addlinespace
 \addlinespace
  Vertices$\backslash$Types     &    All &  $S^3$ & $S^2\!\times\!S^1$ & $S^2\hbox{$\times\hspace{-1.62ex}\_\hspace{-.4ex}\_\hspace{.7ex}$}S^1$ \\ \midrule
 \addlinespace
 \addlinespace
      5 &      1 &      1 &        --  &  -- \\
 \addlinespace
      6 &      2 &      2 &        --  &  -- \\
 \addlinespace
      7 &      5 &      5 &        --  &  -- \\
 \addlinespace
      8 &     39 &     39 &        --  &  -- \\
 \addlinespace
      9 &   1297 &   1296 &        --  &   1 \\
 \addlinespace
     10 & 249015 & 247882 &        518 & 615 \\ \bottomrule
\end{tabular*}
\end{table}

In this paper, the enumeration of $3$-manifolds is continued:
We completely classify triangulated $3$-manifolds with $10$ vertices.
Moreover, we determine the combinatorial automorphism groups
of all triangulations with up to $10$ vertices,
and we test for all $3$-spheres (and all $3$-balls) with up to
$10$ vertices (with up to $9$ vertices) whether they are
constructible, shellable, or vertex-decomposable.
(See \cite{SulankeLutz2006pre} for enumeration results for triangulated
$3$-manifolds with $11$ vertices.)

\section{Enumeration}

We used a backtracking approach, described as mixed-lexicographic enumeration in \cite{Lutz2005pre}, 
to determine all triangulated $3$-manifolds with $10$~vertices:
The vertex-links of a triangulated $3$-manifold with $10$~vertices
are triangulated $2$-spheres with up to $9$ vertices.
Altogether, there are $73$ such $2$-spheres, which are processed 
in decreasing size. As a first vertex-star of a $3$-manifold 
that we are going to build we take the cone over one of the 
respective $2$-spheres and then add further tetrahedra 
(in lexicographic order) as long as this is possible.
If, for example, a triangle of a partial complex that we built
is contained in three tetrahedra, then this violates
the \emph{pseudo-manifold property}, which requires that in a triangulated $3$-mani\-fold
every triangle is contained in \emph{exactly two} tetrahedra.
We backtrack, remove the last tetrahedron that we added,
and try to add to our partial complex the next tetrahedron 
(with respect to the lexicographic order). See \cite{Lutz2005pre}
for further details on the enumeration.

\begin{thm}
There are precisely $249015$ triangulated $3$-manifolds with $10$~vertices:
$247882$ of these are triangulated $3$-spheres, 
$518$ are vertex-minimal triangulations of the sphere product $S^2\!\times\!S^1$,
and $615$ are triangulations of the twisted sphere product $S^2\hbox{$\times\hspace{-1.62ex}\_\hspace{-.4ex}\_\hspace{.7ex}$}S^1$.
\end{thm}
Table~\ref{tbl:ten3d_leq10} gives the total numbers of all triangulations with up to $10$ vertices.
The numbers of $10$-vertex triangulations are listed in detail in Table~\ref{tbl:ten3d_10}. 
All triangulations can be found online at \cite{Lutz_PAGE}.
The topological types were determined with the bistellar flip program
BISTELLAR \cite{Lutz_BISTELLAR}; see \cite{BjoernerLutz2000} for a description. 

\enlargethispage*{5mm}

For a given triangulation, it is a purely combinatorial task
to determine its combinatorial symmetry group.
We computed the respective groups with a program
written in GAP \cite{GAP4}.
\begin{cor}
There are exactly $1$, $1$, $5$, $36$, $408$, and $7443$ triangulated $3$-manifolds
with $5$, $6$, $7$, $8$, $9$, and $10$ vertices, respectively, that
have a non-trivial combinatorial symmetry group. 
\end{cor}
The symmetry groups along with the numbers of combinatorial
types of triangulations that correspond to a particular group are listed in Table~\ref{tbl:3man}.
Altogether, there are $14$ examples that have a vertex-transitive symmetry group; see \cite{KoehlerLutz2005pre}.




\begin{table}
\small\centering
\defaultaddspace=0.15em
\caption{Combinatorial $3$-manifolds with $10$ vertices.}\label{tbl:ten3d_10}
\begin{tabular*}{10cm}{@{\extracolsep{\fill}}c@{\hspace{5mm}}r@{\hspace{5mm}}r@{\hspace{5mm}}r@{\hspace{5mm}}r@{}}
 \addlinespace
 \addlinespace
 \addlinespace
 \addlinespace
 \addlinespace
 \addlinespace
\toprule
 \addlinespace
 \addlinespace
$f$-vector$\backslash$Types     &    All &  $S^3$ & $S^2\!\times\!S^1$ & $S^2\hbox{$\times\hspace{-1.62ex}\_\hspace{-.4ex}\_\hspace{.7ex}$}S^1$ \\ \midrule
 \addlinespace
 \addlinespace
 (10,30,40,20) &     30 &     30 &         -- &  -- \\
 \addlinespace
 (10,31,42,21) &    124 &    124 &         -- &  -- \\
 \addlinespace
 (10,32,44,22) &    385 &    385 &         -- &  -- \\
 \addlinespace
 (10,33,46,23) &    952 &    952 &         -- &  -- \\
 \addlinespace
 (10,34,48,24) &   2142 &   2142 &         -- &  -- \\
 \addlinespace
 (10,35,50,25) &   4340 &   4340 &         -- &  -- \\
 \addlinespace
 (10,36,52,26) &   8106 &   8106 &         -- &  -- \\
 \addlinespace
 (10,37,54,27) &  13853 &  13853 &         -- &  -- \\
 \addlinespace
 (10,38,56,28) &  21702 &  21702 &         -- &  -- \\
 \addlinespace
 (10,39,58,29) &  30526 &  30526 &         -- &  -- \\
 \addlinespace
 (10,40,60,30) &  38575 &  38553 &         10 &  12 \\
 \addlinespace
 (10,41,62,31) &  42581 &  42498 &         37 &  46 \\
 \addlinespace
 (10,42,64,32) &  39526 &  39299 &        110 & 117 \\
 \addlinespace
 (10,43,66,33) &  28439 &  28087 &        162 & 190 \\
 \addlinespace
 (10,44,68,34) &  14057 &  13745 &        145 & 167 \\
 \addlinespace
 (10,45,70,35) &   3677 &   3540 &         54 &  83 \\ \midrule
 \addlinespace
 \addlinespace
 Total:        & 249015 & 247882 &        518 & 615 \\ \bottomrule
\end{tabular*}
\end{table}


\bigskip

All simplicial $3$-spheres with up to $7$ vertices are polytopal. 
However, there are two non-polytopal $3$-spheres with $8$ vertices,
the Gr\"unbaum and Sreedharan sphere \cite{GruenbaumSreedharan1967} 
and the Barnette sphere \cite{Barnette1973c}.
The classification of triangulated $3$-spheres with $9$ vertices into polytopal and
non-polytopal spheres was started by Altshuler and Steinberg
\cite{AltshulerSteinberg1973}, \cite{AltshulerSteinberg1974}, \cite{AltshulerSteinberg1976}
and completed by Altshuler, Bokowski, and Steinberg \cite{AltshulerBokowskiSteinberg1980}
and Engel \cite{Engel1991}.
For neighborly simplicial $3$-spheres with $10$ vertices 
the numbers of polytopal and non-polytopal spheres
were determined by Altshuler \cite{Altshuler1977},
Bokowski and Garms \cite{BokowskiGarms1987},
and Bokowski and Sturmfels~\cite{BokowskiSturmfels1987}. 

\begin{probl}
Classify all simplicial $3$-spheres with $10$ vertices
into polytopal and non-poly\-topal spheres.
\end{probl}

\vfill

\pagebreak

{\small
\defaultaddspace=.1em

\setlength{\LTleft}{0pt}
\setlength{\LTright}{0pt}
\begin{longtable}{@{}l@{\extracolsep{10pt}}l@{\extracolsep{10pt}}r@{\extracolsep{14pt}}l@{\extracolsep{10pt}}r@{\extracolsep{\fill}}l@{\extracolsep{10pt}}l@{\extracolsep{10pt}}r@{\extracolsep{14pt}}l@{\extracolsep{10pt}}r@{}}
\caption{\protect\parbox[t]{15cm}{Symmetry groups of triangulated $3$-manifolds with up to $10$ vertices.}}\label{tbl:3man}
\\\toprule
 \addlinespace
 \addlinespace
 \addlinespace
 \addlinespace
 $n$ &   Manifold   & $|G|$ &   $G$   & Types &     $n$ &   Manifold  &  $|G|$  &   $G$ &  Types   \\ \cmidrule{1-5}\cmidrule{6-10}
\endfirsthead
\caption{\protect\parbox[t]{15cm}{Symmetry groups of triangulated $3$-manifolds (continued).}}
\\\toprule
 \addlinespace
 \addlinespace
 \addlinespace
 \addlinespace
 $n$ &   Manifold   & $|G|$ &   $G$   & Types &     $n$ &   Manifold  &  $|G|$  &   $G$ &  Types   \\ \cmidrule{1-5}\cmidrule{6-10}
\endhead
\bottomrule
\endfoot
 \addlinespace
 \addlinespace
 \addlinespace
  5  & $S^3$                  & 120 & $S_5$,                              &        & 10  & $S^3$                  &   1 & trivial                             & 240683 \\
 \addlinespace
     &                        &     & transitive                          &      1 &     &                        &   2 & ${\mathbb Z}_2$                     &   6675 \\
 \addlinespace
     &                        &     &                                     &        &     &                        &   3 & ${\mathbb Z}_3$                     &     10 \\
 \addlinespace
  6  & $S^3$                  &  48 & $O^*={\mathbb Z}_2\wr S_3$          &      1 &     &                        &   4 & ${\mathbb Z}_4$                     &     53 \\
 \addlinespace
     &                        &  72 & $S_3\wr{\mathbb Z}_2$,              &        &     &                        &     & ${\mathbb Z}_2\times {\mathbb Z}_2$ &    358 \\
 \addlinespace
     &                        &     & transitive                          &      1 &     &                        &   5 & ${\mathbb Z}_5$                     &      1 \\
 \addlinespace
     &                        &     &                                     &        &     &                        &   6 & ${\mathbb Z}_6$                     &      1 \\
 \addlinespace
  7  & $S^3$                  &   8 & $D_4$                               &      2 &     &                        &     & $S_3$                               &     19 \\
 \addlinespace
     &                        &  12 & $S_3\times {\mathbb Z}_2$           &      1 &     &                        &   8 & ${\mathbb Z}_2^{\,3}$               &     15 \\
 \addlinespace
     &                        &  14 & $D_7$,                              &        &     &                        &     & $D_4$                               &     31 \\
 \addlinespace
     &                        &     & transitive                          &      1 &     &                        &  10 & ${\mathbb Z}_{10}$,                 &        \\
 \addlinespace
     &                        &  48 & $D_4\times D_3$                     &      1 &     &                        &     & transitive                          &      1 \\
 \addlinespace
     &                        &     &                                     &        &     &                        &     & $D_5$                               &      4 \\
 \addlinespace
  8  & $S^3$                  &   1 & trivial                             &      3 &     &                        &  12 & $S_3\times {\mathbb Z}_2$           &     15 \\
 \addlinespace
     &                        &   2 & ${\mathbb Z}_2$                     &     13 &     &                        &  16 & $D_4\times {\mathbb Z}_2$           &      3 \\
 \addlinespace
     &                        &   4 & ${\mathbb Z}_4$                     &      1 &     &                        &  20 & $D_{10}$,                           &        \\
 \addlinespace
     &                        &     & ${\mathbb Z}_2\times {\mathbb Z}_2$ &      9 &     &                        &     & transitive                          &      1 \\
 \addlinespace
     &                        &   6 & $S_3$                               &      1 &     &                        &     & $AGL(1,5)$,                         &        \\
 \addlinespace
     &                        &   8 & ${\mathbb Z}_2^{\,3}$               &      1 &     &                        &     & transitive                          &      2 \\
 \addlinespace
     &                        &     & $D_4$                               &      3 &     &                        &  24 & $T^*=S_4$                           &      1 \\
 \addlinespace
     &                        &  12 & $S_3\times {\mathbb Z}_2$           &      4 &     &                        &     & $D_6\times {\mathbb Z}_2$           &      2 \\
 \addlinespace
     &                        &  16 & $D_4\times {\mathbb Z}_2$           &      1 &     &                        &  48 & $O^*={\mathbb Z}_2\wr S_3$          &      2 \\
 \addlinespace
     &                        &     & $D_8$,                              &        &     &                        &  84 & $D_7\times D_3$                     &      1 \\
 \addlinespace
     &                        &     & transitive                          &      1 &     &                        &  96 & $D_6\times D_4$                     &      1 \\
 \addlinespace
     &                        &  60 & $D_5\times D_3$                     &      1 &     &                        & 120 & $S_5$                               &      1 \\
 \addlinespace
     &                        & 384 & ${\mathbb Z}_2\wr S_4$,             &        &     &                        & 200 & $D_5\wr{\mathbb Z}_2$,              &        \\
 \addlinespace
     &                        &     & transitive                          &      1 &     &                        &     & transitive                          &      1 \\
 \addlinespace
     &                        &     &                                     &        &     &                        & 240 & $S_5\times {\mathbb Z}_2$,          &        \\
 \addlinespace
  9  & $S^3$                  &   1 & trivial                             &    889 &     &                        &     & transitive                          &      1 \\
 \addlinespace
     &                        &   2 & ${\mathbb Z}_2$                     &    319 &     & $S^2\!\times\!S^1$     &   1 & trivial                             &    420 \\
 \addlinespace
     &                        &   3 & ${\mathbb Z}_3$                     &      3 &     &                        &   2 & ${\mathbb Z}_2$                     &     95 \\
 \addlinespace
     &                        &   4 & ${\mathbb Z}_4$                     &      6 &     &                        &  10 & ${\mathbb Z}_{10}$,                 &        \\
 \addlinespace
     &                        &     & ${\mathbb Z}_2\times {\mathbb Z}_2$ &     46 &     &                        &     & transitive                          &      1 \\
 \addlinespace
     &                        &   6 & ${\mathbb Z}_6$                     &      1 &     &                        &  16 & $\langle\, 2,2,2\,\rangle_2$        &      1 \\
 \addlinespace
     &                        &     & $S_3$                               &      8 &     &                        &  20 & $D_{10}$,                           &        \\
 \addlinespace
     &                        &   8 & ${\mathbb Z}_2^{\,3}$               &      3 &     &                        &     & transitive                          &      1 \\
 \addlinespace
     &                        &     & $D_4$                               &      5 &     & $S^2\hbox{$\times\hspace{-1.62ex}\_\hspace{-.4ex}\_\hspace{.7ex}$}S^1$  &   1 & trivial                             &    469 \\
 \addlinespace
     &                        &  12 & $S_3\times {\mathbb Z}_2$           &     10 &     &                        &   2 & ${\mathbb Z}_2$                     &    127 \\
 \addlinespace
     &                        &  18 & $D_9$,                              &        &     &                        &   4 & ${\mathbb Z}_2\times {\mathbb Z}_2$ &     14 \\
 \addlinespace
     &                        &     & transitive                          &      1 &     &                        &   8 & $D_4$                               &      2 \\
 \addlinespace
     &                        &  24 & $T^*=S_4$                           &      3 &     &                        &  10 & $D_5$                               &      1 \\
 \addlinespace
     &                        &  72 & $D_6\times D_3$                     &      1 &     &                        &  20 & $D_{10}$,                           &        \\
 \addlinespace
     &                        &  80 & $D_5\times D_4$                     &      1 &     &                        &     & transitive                          &      2 \\
 \addlinespace
     & $S^2\hbox{$\times\hspace{-1.62ex}\_\hspace{-.4ex}\_\hspace{.7ex}$}S^1$ &  18 & $D_9$,                              &        \\
 \addlinespace
     &                        &     & transitive                          &      1 \\
 \addlinespace
 \addlinespace
 \addlinespace
 \addlinespace
\end{longtable}

}

\section{3-Balls}

Along with the enumeration of triangulated $3$-spheres with up to $10$ vertices 
we implicitly enumerated all triangulated $3$-balls with up to $9$ vertices:
Let $B^3_{n-1}$ be a triangulated $3$-ball with $n-1$ vertices and let $v_n$ be a new vertex. 
Then the union\, $B^3_{n-1}\cup (v_n*\partial B^3_{n-1})$\, of $B^3_{n-1}$ with the cone\, $v_n*\partial B^3_{n-1}$\, 
over the boundary $\partial B^3_{n-1}$ with respect to $v_n$ is a triangulated $3$-sphere.
Thus there are at most as many combinatorially distinct $3$-spheres
with $n$ vertices as there are combinatorially distinct $3$-balls with
$n-1$ vertices. 
If, on the contrary, we delete the star of a
vertex from a triangulated $3$-sphere $S^3_n$ with $n$ vertices, 
then, obviously, we obtain a $3$-ball with $n-1$ vertices.
If we delete the star of a different vertex from $S^3_n$ then we
might or might not obtain a combinatorially different ball. 
Let $\#B^3(n-1)$ and $\#S^3(n)$ be the numbers of combinatorially 
distinct $3$-balls and $3$-spheres with $n-1$ and $n$ vertices, respectively.
Then
$$\#S^3(n)\leq \#B^3(n-1)\leq n\cdot \#S^3(n).$$
For the explicit numbers of simplicial $3$-balls with up to $9$ vertices see Table~\ref{tbl:ten3d_balls}.

\begin{table}
\small\centering
\defaultaddspace=0.15em
\caption{Combinatorial $3$-balls with up to $9$ vertices.}\label{tbl:ten3d_balls}
\begin{tabular*}{\linewidth}{@{\extracolsep{\fill}}c@{\hspace{5mm}}r@{\hspace{5mm}}r@{\hspace{5mm}}r@{}}
 \addlinespace
 \addlinespace
 \addlinespace
 \addlinespace
 \addlinespace
 \addlinespace
\toprule
 \addlinespace
 \addlinespace
  Vertices$\backslash$Types     &    All &  Non-Shellable & Not Vertex-Decomposable \\ \midrule
 \addlinespace
 \addlinespace
      4 &       1 &   -- &     -- \\
 \addlinespace
      5 &       3 &   -- &     -- \\
 \addlinespace
      6 &      12 &   -- &     -- \\
 \addlinespace
      7 &     167 &   -- &      2 \\
 \addlinespace
      8 &   10211 &   -- &    628 \\
 \addlinespace
      9 & 2451305 &   29 & 623819 \\ \bottomrule
\end{tabular*}
\end{table}

\section{Vertex-Decomposability, Shellability, \\ and Constructibility}

The concepts of \emph{vertex-decomposability}, \emph{shellability},
and \emph{constructibility} describe three particular ways to assemble
a simplicial complex from the collection of its facets (cf.\ Bj\"orner
\cite{Bjoerner1995} and see the surveys \cite{HachimoriZiegler2000},
\cite{Lutz2004c}, and \cite{Ziegler1998}).
The following implications are strict for (pure) simplicial complexes:
\begin{center}
vertex decomposable\, $\Longrightarrow$\, shellable\, $\Longrightarrow$\, constructible.
\end{center}
It follows from Newman's and Alexander's fundamental works
on the foundations of combinatorial and PL topology from 
1926 \cite{Newman1926c} and 1930 \cite{Alexander1930} 
that a constructible $d$-dimensional simplicial complex 
in which every $(d-1)$-face is contained in exactly two 
or at most two $d$-dimensional facets is a PL $d$-sphere 
or a PL $d$-ball, respectively.

A \emph{shelling} of a triangulated $d$-ball or $d$-sphere 
is a linear ordering of its $f_d$ facets $F_1,...,F_{f_d}$ 
such that if we remove the facets from the ball or sphere in this
order, then at every intermediate step the remaining simplicial
complex is a simplicial ball. A simplicial ball or sphere is \emph{shellable}
if it has a shelling; it is \emph{extendably shellable} 
if any partial shelling $F_1,...,F_i$, $i<f_d$, can be extended to a
shelling; and it is \emph{strongly non-shellable} 
if it has no \emph{free} facet that can be removed from the triangulation without loosing ballness.

A triangulated $d$-ball or $d$-sphere is \emph{constructible}
if it can be decomposed into two constructible $d$-balls of smaller size
(with a single $d$-simplex being constructible)
and if, in addition, the intersection of the two balls is a
constructible ball of dimension $d-1$.
A triangulated $d$-ball or $d$-sphere is \emph{vertex-decomposable} if we can remove the star of 
a vertex $v$ such that the remaining complex is a vertex-decomposable $d$-ball
(with a single $d$-simplex being vertex-decomposable)
and such that the link of $v$ is a vertex-decomposable $(d-1)$-ball or a vertex-decomposable 
$(d-1)$-sphere, respectively.

We tested vertex-decomposability and shellability with
a straightforward backtracking implementation.

\begin{cor}
All triangulated $3$-spheres with $n\leq 10$ vertices are shellable
and therefore constructible.
\end{cor}

An example of a non-constructible and thus non-shellable $3$-sphere
with $13$ vertices was constructed in \cite{Lutz2004c}, whereas
all $3$-spheres with $11$ vertices \emph{are} shellable \cite{SulankeLutz2006pre}.
It remains open whether there are non-shellable respectively
non-constructible $3$-spheres with $12$ vertices.

\begin{cor}
All triangulated $3$-balls with $n\leq 8$ vertices are shellable
and therefore extendably shellable.
\end{cor}
Examples of non-shellable $3$-balls can be found at various places in
the literature (cf.\ the references in \cite{Lutz2004c},
\cite{Lutz2004a}, and \cite{Ziegler1998}) with the smallest previously
known non-shellable $3$-ball by Ziegler \cite{Ziegler1998} with $10$ vertices.

\begin{cor}
There are precisely $29$ vertex-minimal non-shellable 
simplicial $3$-balls with $9$ vertices, ten of which are strongly
non-shellable. The twenty-nine balls have between $18$ and $22$ facets, 
with one unique ball $B\_3\_9\_18$ having $18$ facets and $f$-vector $(9,33,43,18)$. 
\end{cor}
A list of the facets and a visualization of the ball $B\_3\_9\_18$ is given in \cite{Lutz2004a}.

The cone over a simplicial $d$-ball with respect to a new vertex
is a $(d+1)$-dimensional ball.
It is shellable respectively vertex-decomposable if and only
if the original ball is shellable respectively vertex-decomposable (cf.\ \cite{ProvanBillera1980}).

\begin{cor}
There are non-shellable $3$-balls with $d+6$ vertices and $18$ facets for $d\geq 3$.
\end{cor}

Each of the $29$ non-shellable $3$-balls with $9$ vertices
can be split into a pair of shellable balls.

\begin{cor}
All triangulated $3$-balls with $n\leq 9$ vertices are constructible.
\end{cor}

Klee and Kleinschmidt \cite{KleeKleinschmidt1987} showed that all 
simplicial $d$-balls with up to $d+3$ vertices are vertex-decomposable.

\begin{cor}
There are not vertex-decomposable $3$-balls with $d+4$ vertices and $10$ facets for $d\geq 3$.
\end{cor}
In fact, there are exactly two not vertex-decomposable $3$-balls with
$7$ vertices; see \cite{Lutz2004b} for a visualization of these two balls.
One of the examples has $10$ tetrahedra, the other has $11$ tetrahedra.

For the numbers of not vertex-decomposable $3$-balls with up to $9$ vertices
see Table~\ref{tbl:ten3d_balls}.

\begin{cor}
All triangulated $3$-spheres with $n\leq 8$ vertices are vertex-decomposable.
\end{cor}

Klee and Kleinschmidt \cite{KleeKleinschmidt1987} constructed an
example of a not vertex-decomposable polytopal $3$-sphere with $10$
vertices.

\begin{cor}
There are precisely $7$ not vertex-decomposable $3$-spheres with $9$
vertices, which are all non-polytopal. Moreover, there are $14468$ 
not vertex-decomposable $3$-spheres with $10$ vertices.
\end{cor}

Four of the seven examples with $9$ vertices are neighborly with $27$ tetrahedra, 
the other three have $25$, $26$, and $26$ tetrahedra, respectively.
The $25$ tetrahedra of the smallest example are:

{\small
\begin{center}
\begin{tabular}{@{\extracolsep{4mm}}lllll}
$1234$ & $1235$ & $1246$ & $1257$ & $1268$ \\
$1278$ & $1345$ & $1456$ & $1567$ & $1679$ \\
$1689$ & $1789$ & $2348$ & $2359$ & $2378$ \\
$2379$ & $2468$ & $2579$ & $3458$ & $3568$ \\
$3569$ & $3689$ & $3789$ & $4568$ & $5679$. 

\end{tabular}
\end{center}
}

\bibliography{.}

\bigskip
\bigskip
\medskip


\noindent
\normalsize
Frank H. Lutz\\
Technische Universit\"at Berlin\\
Fakult\"at II - Mathematik und Naturwissenschaften\\
Institut f\"ur Mathematik, Sekr.\ MA 3-2\\
Stra\ss e des 17.\ Juni 136\\
10623 Berlin\\
Germany\\
{\tt lutz@math.tu-berlin.de}

\end{document}